\documentclass[letter, 10pt, conference]{ieeeconf}
\IEEEoverridecommandlockouts        

\usepackage{etex}
\usepackage[vlined,ruled]{algorithm2e}
\usepackage[dvipsnames]{xcolor}
\usepackage{pgfpages}
\usepackage[cmex10]{amsmath}		
\usepackage{amssymb, mathrsfs}
\usepackage{cite}
\usepackage{url}

\usepackage{dsfont}
\usepackage{siunitx}
\usepackage{verbatim}									
\usepackage{arydshln}
\usepackage{mathtools}									
\usepackage{siunitx}									
\usepackage{mdframed}

\usepackage{graphicx}
\usepackage{epstopdf}

\usepackage[font=small]{caption}
\usepackage{subcaption}
\usepackage[activate={true,nocompatibility},final,tracking=true,kerning=true,spacing=true,factor=1100,stretch=10,shrink=10]{microtype}

\usepackage{array}
\usepackage{multirow}
\usepackage{enumerate}
\usepackage{booktabs}

\graphicspath{{./fig/}}

\usepackage{pgfplots}
\pgfplotsset{width=8cm,compat=newest}
\newcommand*{\damping}{0.006}%
\newcommand*{\freq}{25}%
\pgfmathsetmacro{\freqd}{sqrt(1-(\damping)^2)*\freq}%

\pgfplotsset{
    standard/.style={
    axis x line=middle,
    axis y line=middle,
    enlarge x limits=0.15,
	enlarge y limits=0.15,
	every axis plot post/.style={mark options={fill=black}},
	}
}

\pgfplotsset{%
    ,compat=1.12
    ,every axis x label/.style={at={(current axis.right of origin)},anchor=north west}
    ,every axis y label/.style={at={(current axis.above origin)},anchor=north east}
    }

\usepackage{tikz}
\usetikzlibrary{arrows}
\usetikzlibrary{decorations.pathmorphing}
\usetikzlibrary{decorations.markings}
\usetikzlibrary{snakes}
\usetikzlibrary{matrix}
\usetikzlibrary{calc}
\usetikzlibrary{patterns}
\usetikzlibrary{positioning}
\usetikzlibrary{shapes}
\usetikzlibrary{shapes.symbols}
\usetikzlibrary{shapes.geometric}
\usetikzlibrary{shadows.blur}
\usetikzlibrary{shapes.arrows}
\usetikzlibrary{shapes.multipart}
\usetikzlibrary{shapes.callouts}
\usetikzlibrary{shapes.misc}

\tikzstyle{every node}=[font=\small]
\tikzstyle{every path}=[line width=0.8pt,line cap=round,line join=round]

\usepackage[american,smartlabels]{circuitikz}
\ctikzset{bipoles/thickness=1}
\ctikzset{bipoles/length=0.8cm}
\ctikzset{bipoles/diode/height=.375}
\ctikzset{bipoles/diode/width=.3}
\ctikzset{tripoles/thyristor/height=.8}
\ctikzset{tripoles/thyristor/width=1}
\ctikzset{bipoles/vsourceam/height/.initial=.7}
\ctikzset{bipoles/vsourceam/width/.initial=.7}

\newcommand{\real}{\mathbb{R}}

\newcommand{\setdef}[2]{\{#1 \;|\; #2\}}

\newcommand{\union}{\operatorname{\cup}}

\DeclareMathOperator*{\minimize}{minimize} 									

\newcommand{\vect}[1]{\mathbbold{#1}}
\newcommand{\vones}[1][]{\vect{1}_{#1}}

\DeclareSymbolFont{bbold}{U}{bbold}{m}{n}
\DeclareSymbolFontAlphabet{\mathbbold}{bbold}

\newcommand{\map}[3]{#1: #2 \rightarrow #3}

\renewcommand{\theenumi}{(\roman{enumi})}

\usepackage{soul}

\newcommand{\define}{\coloneqq}

\DeclareMathOperator{\subto}{subject~to}

\newcommand\oprocendsymbol{\hbox{$\square$}}
\newcommand\oprocend{\relax\ifmmode\else\unskip\hfill\fi\oprocendsymbol}

\newtheorem{theorem}{Theorem}[section]
\newtheorem{lemma}[theorem]{Lemma}

\newtheorem{remark}{Remark}[section]
\newtheorem{assumption}{Assumption}[section]

\newenvironment{pfof}[1]{\vspace{1ex}\noindent{\itshape Proof of
    #1:}\hspace{0.5em}} {\hfill\oprocend\vspace{1ex}}

\newif\ifforstudents

\newcommand{\T}{\mathsf{T}}

\binoppenalty=\maxdimen
\relpenalty=\maxdimen

\title{\bf On Stability of Distributed-Averaging Proportional-Integral Frequency Control in Power Systems
\thanks{
}
}

\author{John W. Simpson-Porco%
    \thanks{
      J. W. Simpson-Porco is with the Department of Electrical and Computer Engineering, University of Waterloo, ON, Canada. Email: {\tt jwsimpson@uwaterloo.ca}. This work was funded under  NSERC Discovery Grant RGPIN-2017-04008 and by UW ECE start-up funding. 
}
}

\begin{document}
\maketitle
\thispagestyle{empty}
\pagestyle{empty}


\begin{abstract}
Distributed consensus-based controllers for optimal secondary frequency regulation of microgrids and power systems have received substantial attention in recent years. This paper provides a Lyapunov-based proof that, under a time-scale separation, these control schemes are stabilizing for a wide class of nonlinear power system models, and under weak assumptions on (i) the objective functions used for resource allocation, and (ii) the graph topology describing communication between agents in the consensus protocol. The results are illustrated via simulation on a detailed test system.
\end{abstract}

\section{Introduction}\label{Sec:Introduction}

The frequency of operation in an AC power system must remain very close to its nominal set-point value of 50 or 60Hz for most equipment to function properly. Due to a combination of the natural physics of synchronous machines, the aggregate behaviour of motors loads, and the conventional \emph{primary control} loops implemented in the system, there is a direct linear relationship between the steady-state frequency deviation present in the system and the mismatch between scheduled generation and demand. The problem of \emph{secondary frequency regulation} is to rebalance  supply and demand, and thereby eliminate any frequency deviation.

The traditional control architecture \cite{PK:94} for achieving secondary frequency regulation is a straight-forward centralized integral control approach: a frequency deviation measurement integrated to produce an overall control signal, which is then allocated to the controllable devices in the system according to so-called participation factors. Presently however, the proliferation of distributed energy resources, flexible loads, and high-bandwidth communication throughout modern small and large-scale power systems has prompted the investigation of alternative distributed control architectures which do not require central coordination. Recent surveys of techniques in this direction include \cite{FD-SB-JWSP-SG:19a,YK-QS-RH-TD-MN-JWSP-FD-MF-FB-HB:19b,DKM-FD-HS-SHL-SC-RB-JL:17}.

This note focuses on one such distributed control scheme, known as the \emph{proportional-integral distributed-averaging} (DAPI) controller. Roughly speaking, this controller uses a multi-agent consensus algorithm to distribute integral action across many controllable devices. The controller was proposed independently in \cite{JWSP-FD-FB:12u,MA-DVD-KHJ-HS:13} as a consensus-based framework for sharing power between generation units and eliminating frequency deviations in microgrids and power systems, respectively. Shortly thereafter, the controller was placed into a distributed optimization framework \cite{FD-JWSP-FB:13y} and experimentally tested for microgrid control \cite{JWSP-QS-FD-JMV-JMG-FB:13s}.

A broad set of literature spanning power electronics, power systems, industrial electronics, and control has subsequently developed around DAPI control; we focus here on the control literature. In terms of stability analysis, \cite{JWSP-FD-FB:12u,MA-DVD-KHJ-HS:13} showed local exponential stability for simple first and second-order power system models, respectively. A Lyapunov-based proof of asymptotic stability for a nonlinear swing-type model appeared in \cite{CZ-EM-FD:15}. A similar proof method can be extended to include basic voltage dynamics \cite{ST-MB-CDP:16}, turbine-governor dynamics \cite{ST-CDP:17,AK-NM-IL:18}, and to show exponential stability \cite{EW-CDP-NM:18}. Other works have analyzed $\mathcal{H}_2$ performance \cite{ET-MA-JWSP-HS:15d,MA-ET-HS-KHJ:17,BKP-JWSP-NM-FD:19d} and performance degradation due to non-cooperative agents \cite{CDP-NM-JWSP:16d}, designed optimal communication topologies \cite{XW-FD-MRJ:16}, examined stability robustness to communication delays \cite{JS-FD-EF:17}, and have placed DAPI within a broad class of optimizing feedback controls \cite{LSPL-JWSP-EM:18l} for linear time-invariant systems.

All of the stability proofs described above are highly dependent on the particular power system model under consideration, as the Lyapunov functions used are constructed to exploit underlying passivity properties of the models. In practice, open-loop AC power system dynamics are stable, but highly uncertain. Not only does this preclude the explicit construction of Lyapunov functions, but in practice, it forces system operators to always use slow, low-gain secondary control schemes. In addition, the above quoted stability results are only applicable to the case of the DAPI controller using (strongly convex) quadratic objective functions and undirected (i.e., bidirectional) communication between agents.



\smallskip

\emph{Contributions:} Our contribution here is to provide a technical proof that the DAPI controller for secondary frequency regulation is asymptotically stabilizing for very general power system models and under substantially weakened assumptions on the objective functions and the inter-agent communication topology. Roughly speaking, the criteria are (i) the power system model need only be asymptotically stable with the steady-state frequency deviation being an affine function of the total injected power, (ii) the objective functions in the DAPI control scheme need only be differentiable and strongly convex (which permits barrier functions), and (iii) the communication graph need only contain a globally reachable node. The key technical insight is that the reduced dynamics obtained after a time-scale separation can be transformed into a nonlinear cascade, which admits a composite-type Lyapunov function. We validate our results via simulation on a detailed 14-machine test system modelling the Australian grid.

\smallskip

\emph{Paper Organization:} Section \ref{Sec:Prelim} records some necessary material on convex functions, graphs, and Laplacian matrices. Section \ref{Sec:PowerSystem} describes the power system model, defines the optimal frequency regulation problem, and presents a preliminary lemma. The main stability result is in Section \ref{Sec:MainResult}. Simulation results on a detailed test system are reported in Section \ref{Sec:Simulations}, with conclusions in Section \ref{Sec:Conclusion}.


\section{Preliminary Material}
\label{Sec:Prelim}

\subsection{Strictly convex functions and their conjugates}
\label{Sec:Convex}

Let $I \subseteq \real$ be a closed interval with non-empty interior, and let $\map{f}{I}{\real}$ be continuously differentiable on $\mathrm{interior}(I)$. We say $f$ is \emph{essentially strictly convex} on $I$ if
\begin{equation}\label{Eq:StrictlyConvex}
(\nabla f(x) - \nabla f(x^{\prime}))(x-x^{\prime}) > 0
\end{equation}
for all $x,x^{\prime} \in \mathrm{interior}(I)$ with $x \neq x^{\prime}$. Note that if $f$ is strictly convex, then $\map{\nabla f}{I}{\real}$ is injective on $\mathrm{interior}(I)$. Again under continuously differentiability, we say $f$ is \emph{essentially smooth} on $I$ if $|f(x_k)| \to +\infty$ whenever $x_k \to x \in \mathrm{boundary}(I)$. The \emph{conjugate} $f^*$ of $f$ is defined as $f^*(p) = \inf_{x \in I}\,[f(x) - p^{\T}x]$. A powerful duality result \cite{RTR:97} is that $f$ is essentially strictly convex on $I$ if and only if $f^*$ is essentially smooth on its domain. As a corollary, if $f$ is both essentially strictly convex  and essentially smooth on $I$, then (i) $\mathrm{dom}(f^*) = \real$, (ii) $f^*$ is essentially strictly convex and essentially smooth on $\real$, and (iii) $(\nabla f)^{-1} = \nabla f^*$. 

A stronger version of this duality \cite{XZ:18} occurs when we consider \emph{strong convexity} and \emph{strong smoothness}. We say $\map{f}{I}{\real}$ is strongly convex with parameter $\mu > 0$ if 
\begin{equation}\label{Eq:StronglyConvex}
(\nabla f(x) - \nabla f(x^{\prime}))(x-x^{\prime}) \geq \mu |x-x^{\prime}|^2
\end{equation}
for all $x,x^{\prime} \in \mathrm{interior}(I)$, and when $\mathrm{dom}(f) = \real$, we say that that $f$ is \emph{strongly smooth} with parameter $L > 0$ if $|\nabla f(x)-\nabla f(x^{\prime})| \leq L|x-x^{\prime}|$ for all $x,x^{\prime} \in \real$. A continuously differentiable mapping $\map{f}{\real}{\real}$ is strongly convex if and only if $\map{f^*}{\real}{\real}$ is strongly smooth.




\subsection{Directed graphs, connectivity, and the Laplacian matrix}
\label{Sec:Graphs}

We will require some elements of graph and algebraic graph theory; see \cite{FB-LNS} for background. A weighted directed graph over $m$ nodes is a triple $\mathcal{G} = (\mathcal{R},\mathcal{E},\mathsf{A})$, where $\mathcal{R}$ satisfying $|\mathcal{R}| = m$ is the set of labels for the nodes, $\mathcal{E} \subseteq \mathcal{R} \times \mathcal{R}$ is the set of directed edges specifying the interconnections between nodes, and $\mathsf{A} \in \real^{m \times m}$ is the adjacency matrix, with elements $a_{ij} \geq 0$ satisfying $a_{ij} > 0$ if and only if $(i,j) \in \mathcal{E}$.
The \emph{Laplacian} matrix $\mathsf{L} \in \real^{m \times m}$ associated with $\mathcal{G}$ is defined element-wise as
\[
\ell_{ij} = \begin{cases}
-a_{ij} &\text{if}\quad i \neq j\\
\sum_{\ell \neq i} a_{i\ell} &\text{if}\quad i = j.
\end{cases}
\]
By construction $\mathsf{L}$ has zero row-sums ($\mathsf{L}\vones[m] = 0$), and hence $0$ is an eigenvalue of $\mathsf{L}$ with right-eigenvector $\vones[m]$. All non-zero eigenvalues of $\mathsf{L}$ have positive real part \cite{FB-LNS}. 

The multiplicity of the $0$ eigenvalue of $\mathsf{L}$ is intimately related to the connections between nodes in $\mathcal{G}$. A \emph{directed path} in $\mathcal{G}$ is an ordered sequence of nodes such that any pair of consecutive nodes in the sequence is a directed edge of $\mathcal{G}$. A node $i \in \mathcal{R}$ is said to be \emph{globally reachable} if for any other node $j \in \mathcal{R} \setminus \{i\}$, there exists a directed path in $\mathcal{G}$ which begins at $j$ and terminates at $i$. 
An elegant result is that $0$ is a simple eigenvalue of $\mathsf{L}$ if and only if $\mathcal{G}$ contains a globally reachable node. In this case, the left-eigenvector $\mathsf{w} \in \real^m$ of $\mathsf{L}$ associated with the simple eigenvalue $0$ has nonnegative elements, and $\mathsf{w}_i > 0$ if and only if node $i \in \mathcal{R}$ is globally reachable.

\section{Power System Model and Optimal Frequency Regulation}
\label{Sec:PowerSystem}

\subsection{Power System Model}

The precise dynamical model of the network will not be of primary concern to us; we will assume a very generic nonlinear power system model of the form
\begin{equation}\label{Eq:NonlinearPowerSystem}
\begin{aligned}
\dot{x}(t) &= f(x(t),u(t),w(t)), \quad x(0) = x_0\\
\Delta \omega(t) &= h(x(t),u(t),w(t))
\end{aligned}
\end{equation}
where $x(t) \in \real^n$ is the vector of states, $u(t) \in \real^m$ is the vector of control inputs, and $w(t) \in \real^{n_w}$ is the vector of (piecewise) constant reference signals, disturbances, and unknown parameters. The model \eqref{Eq:NonlinearPowerSystem} may describe a microgrid or a transmission system, and may have been obtained from a more general differential-algebraic model under appropriate regularity conditions. The controls $u$ will represent power injection set-points for resources participating in secondary frequency regulation; we let $\mathcal{R}$ be an index set for these resources. The disturbance $w$ models set-point changes to other control loops and unmeasured load and generation changes, e.g., from renewable sources. The measurable output $\Delta \omega(t) \in \real^m$ is the vector of frequency deviations at the secondary control resources.

\smallskip


\begin{assumption}[\bf Power System Model]\label{Ass:PowerSystem}
There exist domains $\mathcal{X} \subseteq \real^n$ and $\mathcal{I} \subseteq \real^m \times \real^{n_w}$  such that
\begin{enumerate}[1)]
\item $f$ and $h$ are Lipschitz continuous on $\mathcal{X} \times \mathcal{I}$,
\item there exists a differentiable map $\map{\pi_{x}}{\mathcal{I}}{\mathcal{X}}$ which is Lipschitz continuous on $\mathcal{I}$ and satisfies
\[
0 = f(\pi_{x}(u,w),u,w), \qquad \text{for all}\,\,\, (u,w) \in \mathcal{I};
\]
\item there exist constants $c_1, c_2, c_3, c_4 > 0$ and a function
\[
\map{V_{\rm ps}}{\mathcal{X} \times \mathcal{I}}{\real_{\geq 0}}, \quad (x,(u,w)) \mapsto V_{\rm ps}(x,u,w)
\]
which is continuously differentiable in $(x,u)$ and satisfies
\[
\begin{aligned}
c_1 \|x-\pi_x(u,w)\|_2^2 \leq V_{\rm ps}&(x,u,w) \leq c_2 \|x-\pi_x(u,w)\|_2^2\\
\nabla_{x} V_{\rm ps}(x,u,w)^{\T}f(x,u,w) &\leq -c_3 \|x - \pi_{x}(u,w)\|_2^2\\
\|\nabla_{u} V_{\rm ps}(x,u,w)\|_2 &\leq c_4 \|x-\pi_{x}(u,w)\|_2
\end{aligned}
\]
for all $x \in \mathcal{X}$ and $(u,w) \in \mathcal{I}$;
\item the equilibrium input-to-frequency map $\map{\Delta\bar{\omega}}{\mathcal{I}}{\real^m}$ defined by $\Delta \bar{\omega}(u,w) = h(\pi_{x}(u,w),u,w)$ has the form
\begin{equation}\label{Eq:DefofPiSpecialized}
\Delta \bar{\omega}(u,w) = \frac{1}{\beta}\vones[m](\vones[m]^{\T}\bar{u} - d),
\end{equation}
where $\beta > 0$ and $d \in \real$ is the (constant) unmeasured net load disturbance.
\end{enumerate}
\end{assumption}

Assumption (2) above says that associated to each constant control/disturbance pair $(u,w) \in \mathcal{I}$ is a unique (at least, on the set $\mathcal{X}$) equilibrium state $\pi_x(u,w)$. Assumption (3) is a Lyapunov function establishing exponential stability of $\pi_x(u,w) \in \mathcal{X}$. Assumptions (1)--(3) are placed to ensure we can pursue a singular perturbation framework for stability analysis; variations and relaxations are possible. Assumption (4) specifies that that the network achieves \emph{frequency synchronization} in steady-state, with frequency deviations being equal at all nodes in the system. The steady-state value of the common frequency is determined by $\vones[m]^{\T}u - d$, the mismatch between generation and demand.

A very simple model which satisfies Assumption \ref{Ass:PowerSystem} is
\[
\begin{aligned}
\Delta\dot{\theta}_i &= \Delta\omega_i\,,\\
M_i\Delta\dot{\omega}_i &= - \sum_{j=1}^{n} T_{ij}(\Delta \theta_i-\Delta \theta_j) - D_i\Delta\omega_i + \Delta P_{\mathrm{m},i} - d_i \\
T_i \Delta \dot{P}_{\mathrm{m},i} &= -\Delta P_{\mathrm{m},i} - R_{\mathrm{d},i}^{-1}\Delta\omega_i + u_i.
\end{aligned}
\]
for $i \in \{1,\ldots,n\}$ with $\Delta\theta_1 \equiv 0$, which describes a linearized network-reduced model of synchronous machines with first-order turbine governor models. While we refer the reader to \cite{FD-SB-JWSP-SG:19a} and the references therin for details on these kinds of models, we note that for this particular model, the constant $\beta$ in \eqref{Eq:DefofPiSpecialized} is given by $\beta = \sum_{i=1}^{m} D_i + R_{\mathrm{d},i}^{-1}$.

\subsection{Optimal and Distributed Frequency Regulation}
\label{Sec:OptimalFrequencyRegulation}

For the goal of secondary frequency regulation, a typical power system is highly over-actuated, and the operator has flexibility in allocating control actions across many actuators. The desired set-points can be specified via the minimization
\begin{subequations}\label{Eq:OSSFreq}
\begin{align}
\label{Eq:OSSFreq-1}
\minimize_{\bar{u} \in \real^m} &\quad J(\bar{u}) \define \sum_{i\in\mathcal{R}}\nolimits J_{i}(\bar{u}_i)
\\
\label{Eq:OSSFreq-2}
\subto &\quad  0 = \vones[m]^{\T}u - d
\end{align}
\end{subequations}
where $\map{J_{i}}{\mathcal{U}_i}{\real}$ models the disutility of the $i$th secondary power provider, and includes a (smooth) barrier function for enforcing inequality constraints $\bar{u}_i \in \mathcal{U}_i = (\underline{u}_i,\overline{u}_i)$, where $-\infty \leq \underline{u}_i < \overline{u}_i \leq +\infty$. In other words, any limit constraints are directly included in the domain of the function $J_i$. The constraint \eqref{Eq:OSSFreq-2} enforces balance of secondary power injections $\vones[m]^{\T}u$ and unmeasured demand $d$, and by \eqref{Eq:DefofPiSpecialized}, enforces that the steady-state network frequency deviation should be zero. We assume that \eqref{Eq:OSSFreq} is strictly feasible.

\smallskip

\begin{assumption}[\bf Regularity of Objective Functions]\label{Ass:Cost}
Each function $\map{J_i}{\mathcal{U}_i}{\real_{\geq 0}}$ is continuously differentiable, strongly convex on $\mathcal{U}_i$ with parameter $\mu_i > 0$, and satisfies the barrier function properties
\[
\lim_{\xi \searrow \underline{u}_i} J_i(\xi_i) = +\infty, \quad \lim_{\xi \nearrow \overline{u}_i} J_i(\xi_i) = +\infty.
\]
\end{assumption}

\medskip

It follows from Assumption \ref{Ass:Cost} that $J_i$ is essentially strictly convex and essentially smooth on $\mathcal{U}_i$ (Section \ref{Sec:Convex}). The control problem of interest is to design a (distributed) feedback controller which drives the system frequency deviation towards zero while simultaneously ensuring the control inputs converge towards the (unique) primal optimizer of \eqref{Eq:OSSFreq}. The distributed-averaging proportional-integral (DAPI) control scheme combines integral control on local frequency measurements with peer-to-peer communication between secondary control resources to solve this problem. The following preliminary result characterizes the optimal solution of \eqref{Eq:OSSFreq}.

\medskip

\begin{lemma}[\bf Distributed Optimality Conditions]\label{Lem:EquivalentOptimality}
Consider the optimization problem \eqref{Eq:OSSFreq}. Let $\bar{u} \in \real^m$, let $\Delta \bar{\omega}$ be as in \eqref{Eq:DefofPiSpecialized}, and let $\mathcal{G} = (\mathcal{R},\mathcal{E},\mathsf{A})$ be a weighted directed graph with associated Laplacian matrix $\mathsf{L}$. Assume that $\mathcal{G}$ contains a globally reachable node, and let $\mathsf{w} \in \real^{m}_{\geq 0}$ be the left-eigenvector of $\mathsf{L}$ corresponding to its simple eigenvalue at $0$. If $K \succeq 0$ is diagonal matrix such that $\mathsf{w}^{\T}K\vones[m] > 0$, then the following statements are equivalent:
\begin{enumerate}
\item \label{Lem:ReducedError-1} $\bar{u}$ is the unique primal optimizer of \eqref{Eq:OSSFreq};
\item \label{Lem:ReducedError-2} there exists a unique vector $\bar{\eta} \in \mathrm{span}(\vones[m])$ such that
\begin{subequations}\label{Eq:ReformulatedOpt2}
\begin{align}
\label{Eq:ReformulatedOpt2-1}
0 &= K\Delta \bar{\omega}(\bar{u},w) + \mathsf{L}\bar{\eta}\\
\label{Eq:ReformulatedOpt2-2}
\bar{u} &= \nabla J^*(\bar{\eta}),
\end{align}
\end{subequations}
where $J^*(\eta) = \sum_{i\in\mathcal{R}}J_i^*(\eta_i)$ is the conjugate of $J$.
\end{enumerate}
\end{lemma}

\medskip

\begin{proof}
First note that since \eqref{Eq:OSSFreq} is strictly feasible, $J(\bar{u})$ is strongly convex, and the constraint matrix $\vones[n]^{\T}$ in \eqref{Eq:OSSFreq-2} has full row rank, the problem \eqref{Eq:OSSFreq} has a unique primal-dual optimal solution $(\bar{u},\lambda)$ for some $\lambda \in \real$, which satisfies the KKT conditions \eqref{Eq:OSSFreq-2} and 
\begin{equation}\label{Eq:Stationarity}
\nabla J(\bar{u}) = \lambda\vones[m] \quad \Longleftrightarrow \quad \bar{u} = \nabla J^*(\lambda\vones[m]).
\end{equation}
Since $0$ is a simple eigenvalue of $\mathsf{L}$ with right-eigenvector $\vones[m]$, there exists a unique value $\lambda$ satisfying \eqref{Eq:Stationarity} if and only if there exists a unique vector $\bar{\eta} \in \mathrm{span}(\vones[m])$ such that
\begin{subequations}\label{Eq:ReformulatedOpt1}
\begin{align}
\label{Eq:ReformulatedOpt1-1}
0 &= \mathsf{L}\bar{\eta}\\
\label{Eq:ReformulatedOpt1-2}
\bar{u} &= \nabla J^*(\bar{\eta}).
\end{align}
\end{subequations}
In addition, trivially, the constraint \eqref{Eq:OSSFreq-2} holds if and only if
\begin{equation}\label{Eq:ExpandedEquality}
0 = \vones[m]\beta^{-1}(\vones[m]^{\T}\bar{u}-d) = \Delta \bar{\omega}(\bar{u},w),
\end{equation}
where we have used \eqref{Eq:DefofPiSpecialized}. We now claim that \eqref{Eq:ReformulatedOpt1-1} and \eqref{Eq:ExpandedEquality} hold if and only if \eqref{Eq:ReformulatedOpt2-1} holds. That \eqref{Eq:ReformulatedOpt1-1} and \eqref{Eq:ExpandedEquality} imply \eqref{Eq:ReformulatedOpt2-1} is trivial. For the other direction, left-multiply \eqref{Eq:ReformulatedOpt2-1} by $\mathsf{w}^{\T}$ to find that
\[
0 = \mathsf{w}^{\T}(\Delta \bar{\omega}(\bar{u},w) + \mathsf{L}\bar{\eta}) = \mathsf{w}^{\T}K\vones[m]\beta^{-1}(\vones[m]^{\T}\bar{u}-d).
\]
The vector $\mathsf{w}$ is non-negative, and is non-zero since the graph $\mathcal{G}$ has a globally reachable node (Section \ref{Sec:Graphs}), and by assumption $\mathsf{w}^{\T}K\vones[m] \neq 0$. We conclude that $\vones[m]^{\T}\bar{u}-d = 0$, and therefore \eqref{Eq:ExpandedEquality} holds. Substituting this into  \eqref{Eq:ReformulatedOpt2-1}, it follows that \eqref{Eq:ReformulatedOpt1-1} holds, which completes the proof.
\end{proof}

\medskip

Lemma \ref{Lem:EquivalentOptimality} leads naturally to the DAPI controller
\begin{equation}\label{Eq:DAPIVector}
\tau \dot{\eta}(t) = -\Delta \omega(t) - \mathsf{L}\eta(t), \quad u(t) = \nabla J^*(\eta(t)),
\end{equation}
where $\tau > 0$ is a tuning gain. The vector $\eta(t) \in \real^m$ is now the dynamic controller state, and the steady-state frequency vector $\Delta \bar{\omega}$ has been replaced by the real-time frequency measurement vector $\Delta \omega(t)$. In components, \eqref{Eq:DAPIVector} is
\begin{subequations}\label{Eq:DAPI}
\begin{align}
\tau \dot{\eta}_i(t) &= -\Delta \omega_i(t) - \sum_{j=1}^{m}\nolimits a_{ij}(\eta_i(t)-\eta_j(t))\\
u_i(t) &= \nabla J_i^*(\eta_i(t)),
\end{align}
\end{subequations}
which emphasizes that \eqref{Eq:DAPIVector} is a distributed controller.

\smallskip

\begin{remark}[\bf Generalized DAPI Controllers]\label{Rem:GenDAPI}
Lemma \ref{Lem:EquivalentOptimality} strongly suggests that one could insert a gain matrix $K$ in front of $\Delta \omega(t)$ in \eqref{Eq:DAPIVector}. Indeed, intuitively, it does not seem necessary that all nodes explicitly integrate their local frequency measurements. For technical reasons though, our analysis is only applicable to \eqref{Eq:DAPIVector}. The situation where not all agents take local frequency measurements seems theoretically interesting, but does not seem especially important in practice. \hfill \oprocend
\end{remark}


\section{Main Result: Closed-Loop Asymptotic Stability with DAPI Control}
\label{Sec:MainResult}

We now state and prove our main result, that the distributed controller \eqref{Eq:DAPI} leads to stable and optimal frequency regulation of the power system \eqref{Eq:NonlinearPowerSystem}. 

\smallskip

\begin{theorem}[\bf Low-Gain Stability with DAPI Control]\label{Thm:DAPIStable}
Consider the power system model \eqref{Eq:NonlinearPowerSystem} under Assumption \ref{Ass:PowerSystem}, interconnected with the  the DAPI controller \eqref{Eq:DAPI} under Assumption \ref{Ass:Cost}. If the communication graph $\mathcal{G}$ contains a globally reachable node, then there exists $\tau^{\star} > 0$ such that for all $\tau \geq \tau^{\star}$, the unique equilibrium point $(\bar{x},\bar{\eta}) \in \real^{n} \times \real^{m}$ of the closed-loop system is asymptotically stable and $\bar{u} = \nabla J^{*}(\bar{\eta})$ is the unique global primal optimizer of \eqref{Eq:OSSFreq}.
\end{theorem}

\begin{pfof}{Theorem \ref{Thm:DAPIStable}}
For the closed-loop system \eqref{Eq:NonlinearPowerSystem} and \eqref{Eq:DAPIVector}, define the new time variable $\ell = t/\tau$, which leads to the singularly perturbed system
\[
\begin{aligned}
\varepsilon \frac{\mathrm{d}x}{\mathrm{d}\ell} &= f(x,u,w), \qquad
\Delta \omega = h(x,u,w)\\
\frac{\mathrm{d}\eta}{\mathrm{d}\ell} &= -\Delta \omega - \mathsf{L}\eta, \qquad
u = \nabla J^*(\eta),
\end{aligned}
\]
where $\varepsilon = 1/\tau$. We will apply \cite[Theorem 11.3]{HKK:02}, which constructs a quadratic-type Lyapunov function for the interconnection. Due to Assumption \ref{Ass:PowerSystem}, the conditions on the boundary layer system in \cite[Theorem 11.3]{HKK:02} are satisfied, and we must study the reduced dynamics
\begin{equation}\label{Eq:ReducedDynamicsOriginal}
\begin{aligned}
\dot{\eta} &= -\Delta \bar{\omega}(u,w) -\mathsf{L}\eta\\
u &= \nabla J^*(\eta),
\end{aligned}
\end{equation}
where $\Delta \bar{\omega}$ is as given in \eqref{Eq:DefofPiSpecialized} and where $\dot{\eta}$ denotes differentiation with respect to the new temporal variable $\ell$. By Lemma \ref{Lem:EquivalentOptimality} with $K = I_m$, the system  \eqref{Eq:ReducedDynamicsOriginal} possesses a unique equilibrium point $\bar{\eta} \in \mathrm{span}(\vones[m])$. Eliminating $u$ from \eqref{Eq:ReducedDynamicsOriginal}, the dynamics are equivalently given by\footnote{See also \cite{AC-JC:14} for closely related dynamics.}
\begin{equation}\label{Eq:ReducedDynamics}
\dot{\eta} = -\tfrac{1}{\beta}\vones[m]\vones[m]^{\T}\nabla J^*(\eta) - \mathsf{L}\eta + \tfrac{1}{\beta}\vones[m]d.
\end{equation}
Define the nonsingular transformation matrix 
\[
T = \begin{bmatrix}
\vones[m] & V_{\perp}
\end{bmatrix},
\]
where $V_{\perp} \in \real^{m \times (m-1)}$ has columns which form an orthonormal basis for the subspace $\setdef{\eta \in \real^m}{\vones[m]^{\T}\eta = 0}$. Consider the change of state variable
\[
\eta = T\begin{bmatrix}z\\
\delta
\end{bmatrix} = \vones[m]z + V_{\perp}\delta, \quad z \in \real,\,\,\delta \in \real^{m-1}.
\]
It is straightforward to see that
\[
z = \tfrac{1}{m}\vones[m]^{\T}\eta, \qquad \delta = V_{\perp}^{\T}\eta,
\]
and by construction, the unique equilibrium $(\bar{z},\bar{\delta}) = T^{\T}\bar{\eta} =  (\tfrac{1}{m}\vones[m]^{\T}\bar{\eta},0)$ satisfies
\begin{equation}\label{Eq:Equilibrium}
0 = -\tfrac{1}{\beta}\vones[m]^{\T}\nabla J^*(\vones[m]\bar{z}) + \tfrac{1}{\beta}d.
\end{equation}
In the $(z,\delta)$ coordinates, the dynamics \eqref{Eq:ReducedDynamics} become
\begin{equation}\label{Eq:CascadeDynamics}
\begin{aligned}
\dot{z} &= -\tfrac{1}{\beta}\vones[m]^{\T}\nabla J^*(\vones[m]z+ V_{\perp}\delta)  - \tfrac{1}{m}\vones[m]^{\T}\mathsf{L}V_{\perp}\delta + \tfrac{1}{\beta}d\\
\dot{\delta} &= -V_{\perp}^{\T}\mathsf{L}V_{\perp}\delta.
\end{aligned}
\end{equation}
We reformulate the $z$-dynamics in \eqref{Eq:CascadeDynamics} by adding and subtracting the term $\tfrac{1}{\beta}\vones[m]^{\T}\nabla J^*(\vones[m]z)$ and using the equilibrium equation \eqref{Eq:Equilibrium}, which yields the equivalent dynamic model
\begin{equation}\label{Eq:CascadeNonlinear}
\begin{aligned}
\dot{z} &= f_1(z,\delta) = \varphi(z) + \psi(z,\delta)\\
\dot{\delta} &= f_2(\delta)
\end{aligned}
\end{equation}
where
\[
\begin{aligned}
\varphi(z) &\define -\tfrac{1}{\beta}\vones[m]^{\T}\left[\nabla J^*(\vones[m]z)-\nabla J^*(\vones[m]\bar{z})\right]\\
\psi(z,\delta) &\define -\tfrac{1}{\beta}\vones[m]^{\T}\left[\nabla J^*(\vones[m]z + V_{\perp}\delta) -\nabla J^*(\vones[m]z)\right]\\
&\qquad - \tfrac{1}{m}\vones[m]^{\T}\mathsf{L}V_{\perp}\delta\\
f_2(\delta) &\define -V_{\perp}^{\T}\mathsf{L}V_{\perp}\delta.
\end{aligned}
\]
Note that $\psi(z,0) = 0$ for all $z \in \real$. Moreover, since each function $J_i$ is strongly convex with parameter $\mu_i > 0$, we have from Section \ref{Sec:Convex} that $J_i^*$ is strongly smooth with parameter $1/\mu_i$, which implies that
\[
|\psi(z,\delta)| \leq (\underbrace{\tfrac{\sqrt{m}}{\beta}\tfrac{1}{\mu_{\rm min}}\|V_{\perp}\|_2 + \tfrac{1}{\sqrt{m}}\|\mathsf{L}V_{\perp}\|_2}_{\define \kappa})\|\delta\|_2
\]
where $\mu_{\rm min} \define \min_{i \in \mathcal{R}}\mu_i$. The dynamics \eqref{Eq:CascadeNonlinear} are in the form of a \emph{nonlinear cascade}; we will construct a composite Lyapunov function for the cascade.

First consider the driving system $f_2$. Since the graph $\mathcal{G}$ contains a globally reachable node, $\mathsf{L}$ has a simple eigenvalue at $0$ with all other eigenvalues having positive real part (Section \ref{Sec:Graphs}). Note that since $\mathsf{L}\vones[m] = 0$, we have
\[
\begin{aligned}
T^{-1}\mathsf{L}T &= \begin{bmatrix}
\tfrac{1}{\sqrt{m}}\vones[m]^{\T} \\ V_{\perp}^{\T}
\end{bmatrix}\mathsf{L}\begin{bmatrix}
\tfrac{1}{\sqrt{m}}\vones[m] & V_{\perp}
\end{bmatrix}\\
&= \begin{bmatrix}
0 & \tfrac{1}{\sqrt{m}}\vones[m]^{\T}\mathsf{L}V_{\perp}\\
0 & V_{\perp}^{\T}\mathsf{L}V_{\perp}
\end{bmatrix}.
\end{aligned}
\]
It follows that $\mathrm{eig}(\mathsf{L}) = \{0\} \union \mathrm{eig}(V_{\perp}^{\T}\mathsf{L}V_{\perp})$, and it must therefore be that all eigenvalues of $-V_{\perp}^{\T}\mathsf{L}V_{\perp}$ have negative real part. By linear Lyapunov theory, there exists $\rho > 0$ and $P \succ 0$ such that with $W(\delta) = \delta^{\T}P\delta$, we satisfy the dissipation inequality
\begin{equation}\label{Eq:DeltaLyapunov}
\nabla W(\delta)^{\T}f_2(\delta) \leq -\rho \|\delta\|_2^2, \qquad \delta \in \real^{m-1}.
\end{equation}
For the driven system $f_1$, consider the Lyapunov candidate $\map{V}{\real}{\real_{\geq 0}}$ defined as
\begin{equation}\label{Eq:V}
\begin{aligned}
V(z) &= \sum_{i=1}^{m}\nolimits \int_{\bar{z}}^{z} [\nabla J_i(\xi)-\nabla J_i(\bar{z})]\,\mathrm{d}\xi\\
&= \sum_{i=1}^{m}\nolimits \left[J_i^*(z) - J_i^*(\bar{z}) - \nabla J_i^*(\bar{z})(z-\bar{z})\right].
\end{aligned}
\end{equation}
By Assumption \ref{Ass:Cost} each function $J_i$ is essentially strictly convex (in fact, strongly convex) and is essentially smooth. It follows (Section \ref{Sec:Convex}) that $J_i^*$ is also essentially strictly convex and essentially smooth. As each summand in $V(z)$ is the difference between $J_i^*$ and its linear approximation at $\bar{z}$, and it follows that $V$ is positive-definite with respect to $\bar{z}$ \cite[Lemma A.2]{JWSP:17e}. Additionally, by essential smoothness of $J_i^*$, we conclude that $V(z) \to +\infty$ as $|z| \to \infty$, so $V$ is radially unbounded. An easy computation shows that
\[
\nabla V(z) = \vones[m]^{\T}\left[\nabla J^*(\vones[m]z)-\nabla J^*(\vones[m]\bar{z})\right].
\]
Again from Section \ref{Sec:Convex}, since each $J_i^*$ is strictly convex, we conclude that $\nabla V(z) = 0$ if and only if $z = \bar{z}$. For $\alpha > 0$, consider now the composite Lyapunov candidate
\[
\mathcal{V}(z,\delta) = V(z) + \alpha W(\delta),
\]
which is positive definite with respect to $(\bar{z},0)$ and is radially unbounded. Easy calculations now show that
\[
\begin{aligned}
\nabla V(z)^{\T}\phi(z) &= -\tfrac{1}{\beta}|\nabla V(z)|^2\\
\nabla V(z)^{\T}\psi(z,\delta) &\leq \kappa |\nabla V(z)|\|\delta\|_2
\end{aligned}
\]
for all $z \in \real$ and $\delta \in \real^{m-1}$. Combining these with \eqref{Eq:DeltaLyapunov}, we find that along trajectories of \eqref{Eq:CascadeNonlinear}
\[
\dot{\mathcal{V}}(z,\delta) \leq -\begin{bmatrix}
|\nabla V(z)|\\
\|\delta\|_2
\end{bmatrix}^{\T}\begin{bmatrix}
\tfrac{1}{\beta} & -\kappa/2\\
-\kappa/2 & \alpha\rho
\end{bmatrix}
\begin{bmatrix}
|\nabla V(z)|\\
\|\delta\|_2
\end{bmatrix}.
\]
Selecting $\alpha > \beta \kappa^2/(4\rho)$, the right-hand side becomes a negative definite with respect to the equilibrium $(\bar{z},0)$. We conclude that the equilibrium $(\bar{z},0)$ of \eqref{Eq:CascadeNonlinear} \textemdash{} or equivalently, the equilibrium $\bar{\eta}$ of the reduced dynamics \eqref{Eq:ReducedDynamics} \textemdash{} is globally asymptotically stable. All conditions of \cite[Theorem 11.3]{HKK:02} are now satisfied, which completes the proof.
\end{pfof}

An interesting aspect of Theorem \ref{Thm:DAPIStable} is that it imposes the weakest possible time-invariant connectivity assumption one can place on $\mathcal{G}$ to ensure consensus, namely the existence of a globally reachable node \cite{FB-LNS}. Indeed, if $\mathcal{G}$ does not contain a globally reachable node, then $\mathsf{L}$ has at least two eigenvalues at $0$, and the reduced dynamics \eqref{Eq:ReducedDynamics} contain a marginally stable mode. Strongly connected, weight-balanced, and undirected communication graphs are all covered as special cases. The analysis also naturally points to where modifications of the assumptions can be made. For instance, if each $J_i$ is assumed to be defined on all of $\real$ and is assumed to be strongly smooth, then (Section \ref{Sec:Convex}) $J_i^*$ will be strongly convex, and one will be able to use \eqref{Eq:StronglyConvex} to conclude that the reduced dynamics \eqref{Eq:ReducedDynamics} are globally \emph{exponentially} stable (in this case, the simple Lyapunov candidate $V(z) = \tfrac{1}{2}\|z-\bar{z}\|_2^2$ can be used in place of \eqref{Eq:V}).

\section{Simulation on Australian Test System}
\label{Sec:Simulations}

We illustrate our result by simulating the controller \eqref{Eq:DAPI} on a highly detailed dynamic power system model based on the south eastern Austalian system \cite{AM-IK-PB-GS:15}. The model contains 14 synchronous generators, with full-order turbine-governor, excitation, and PSS models. We will use 5 of these generators (buses 201, 301, 401, 403, and 503) as controllable for secondary frequency regulation, with the inputs $u_i$ being the power set-points to their turbine-governor systems.

To expoloit the full flexibility of the theoretical result, we consider heterogeneous objective functions\footnote{The computation of $u_i(t)$ in \eqref{Eq:DAPI} is done by solving the algebraic constraint $\nabla J_i(u_i(t)) = \eta_i(t)$.}
\[
J_i(u_i) =  \tfrac{1}{2}q_i (u_i-u_i^{\star})^2 - \gamma [\log(\overline{u}_i-u_i) +\log(-\underline{u}_i+u_i)]
\]
where $u_i^{\star}$ is the base dispatch point of the resource, $q_i > 0$, and $\gamma > 0$ is a barrier function parameter; the parameters are listed in Table \ref{Tab:Para}. The upper and lower power limits for each resource were chosen as $\pm 0.1$ p.u. from the respective dispatch point. The communication graph $\mathcal{G}$ is an directed line graph (with weights $a_{ij} = 0.1$) connecting the five controllable machines; bus 503 is therefore the unique globally reachable.

\begin{table}[ht!]
\begin{center}
\begin{tabular}{c|ccccc}
\toprule
 & G201 & G301 & G401 & G403 & G503\\
\midrule
$q_i$ & 1 & 0.8 & 1 & 0.8 & 0.1\\
$u_i^{\star}$ (p.u.) & 0.9 & 0.9 & 0.787 & 0.787 & 0.6539\\
\bottomrule
\end{tabular}
\end{center}
\caption{Parameters for simulation study; $\gamma = 0.001$, $\tau = 0.2$.}
\label{Tab:Para}
\end{table}

\begin{figure}[ht!]
\centering
\begin{subfigure}{0.99\linewidth}
\includegraphics[width=\linewidth]{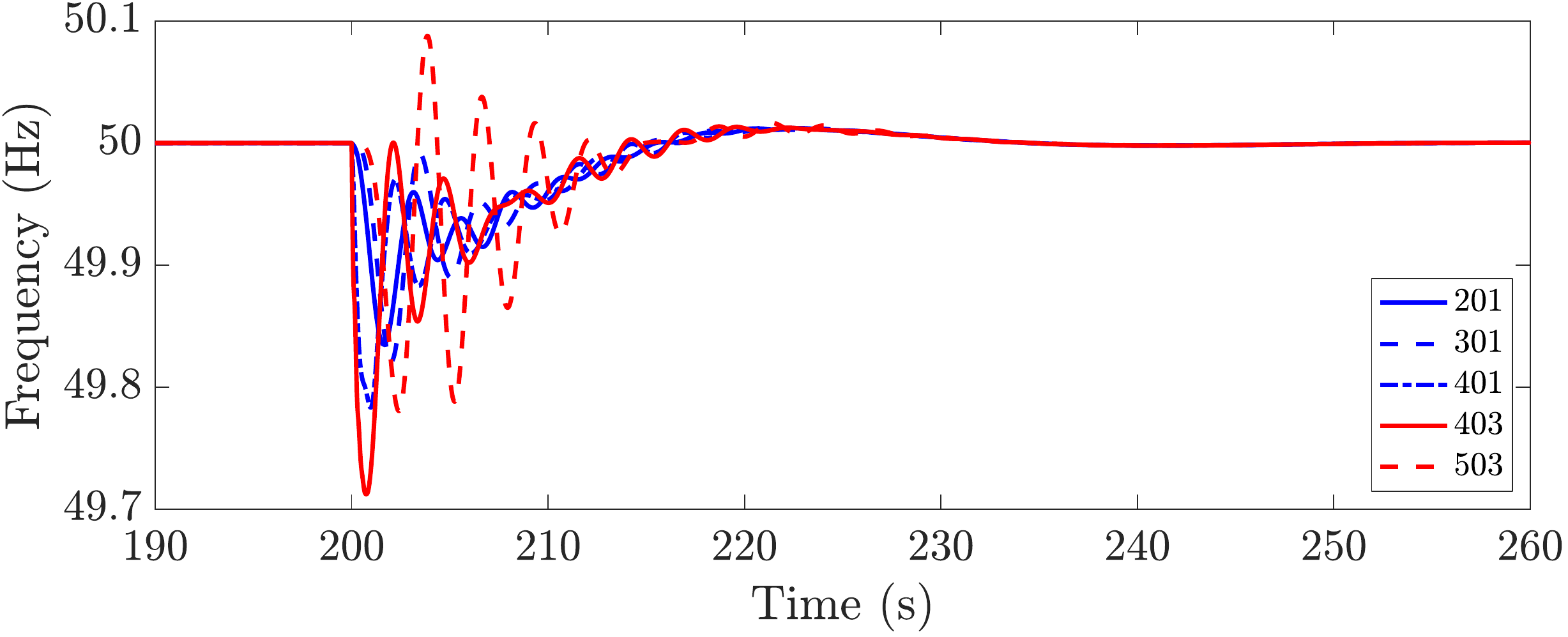}
\caption{Frequency measurements from controllable resources.}
\label{Fig:Aus-a}
\end{subfigure}\\
\smallskip
\begin{subfigure}{0.99\linewidth}
\includegraphics[width=\linewidth]{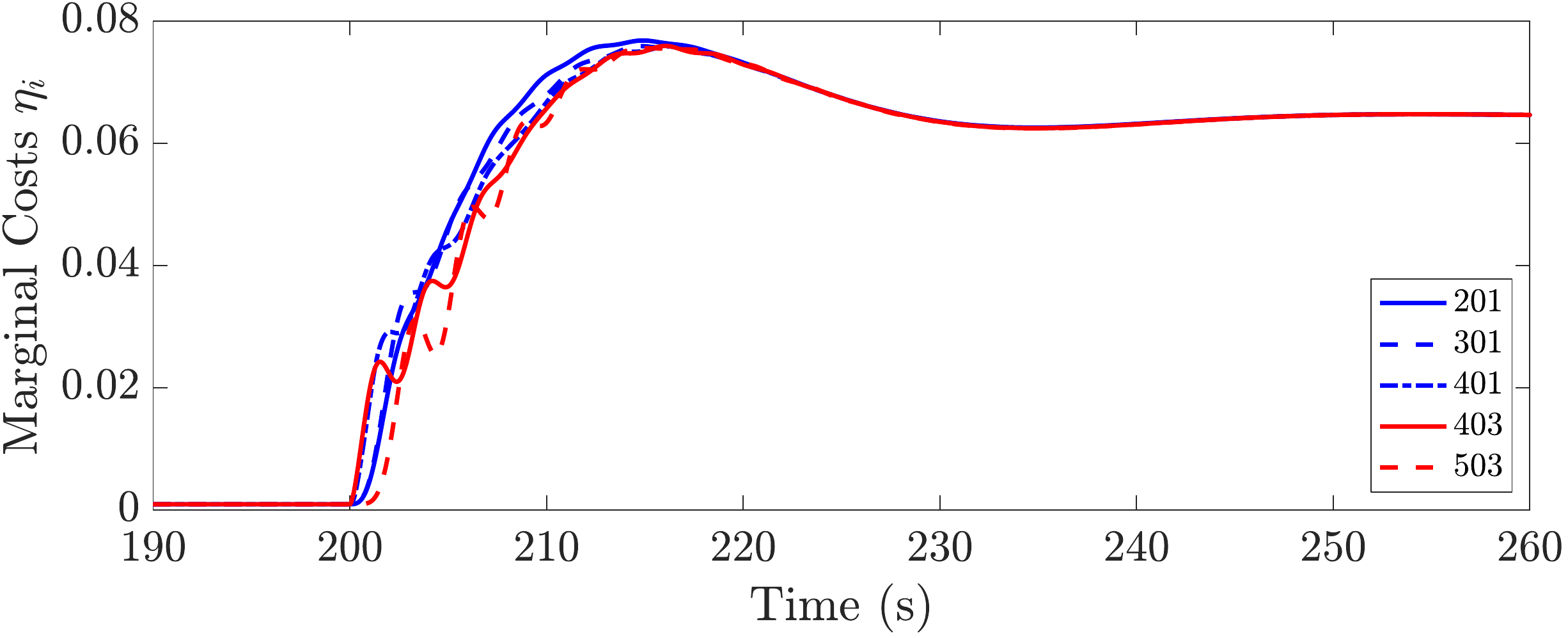}
\caption{Marginal cost variables from \eqref{Eq:DAPI}.}
\label{Fig:Aus-b}
\end{subfigure}\\
\smallskip
\begin{subfigure}{0.99\linewidth}
\includegraphics[width=\linewidth]{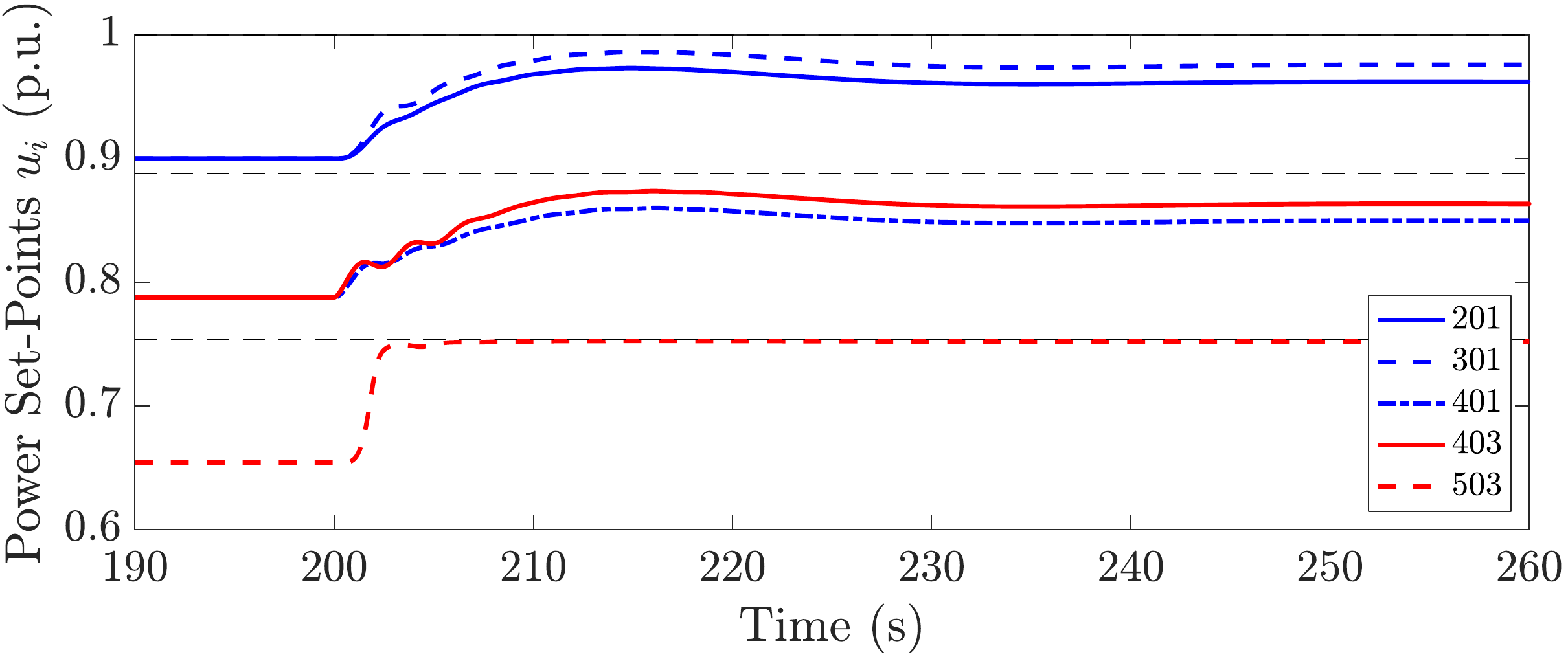}
\caption{Set-points $u_i$ for turbine/governors; black dashed lines are upper limits $\overline{u}_i$.}
\label{Fig:Aus-c}
\end{subfigure}
\caption{Australian 14-machine system with DAPI control.}
\label{Fig:Aus}
\end{figure}

Figure \ref{Fig:Aus} shows the closed-loop response when the load at bus 406 is doubled at time $t = 200$. This is a sizeable disturbance, and inter-area modes are visible in the frequency plot of Figure \ref{Fig:Aus-a}; as expected though, the frequency deviation is asymptotically eliminated. The consensus action in \eqref{Eq:DAPI} keeps the marginal cost variable $\eta_i$ in agreement, as shown in Figure \ref{Fig:Aus-b}. Figure \ref{Fig:Aus-c} shows the set-points $u_i$ sent to the resources. As G503 has a low cost $q_i$ parameter, it is preferentially used, but the log barrier functions ensure that the commands $u_i$ always satisfy the inequality constraints $u_i \in [\underline{u}_i,\overline{u}_i]$.

\section{Conclusions}
\label{Sec:Conclusion}



We have presented a time-scale separation proof that distributed-averaging proportional-integral (DAPI) optimal frequency control is stabilizing for general power system models and under weak assumptions on the objective functions and the inter-agent communication topology uses for consensus. This result closes a persistent gap in the literature and provides credibility for the safe application of \eqref{Eq:DAPI} to practical microgrids and power systems. One remaining open question was noted in Remark \ref{Rem:GenDAPI}. Another unresolved question is how far the convexity assumptions on the objective functions $J_i$ can be relaxed.

\section{Acknowledgements}

This article is dedicated to the memory of Martin Andreasson. The author also acknowledges F. D\"{o}rfler, A. Cherukuri, S. Trip, T. Stegink, E. Tegling, N. Monshizadeh, and C. De Persis for stimulating conversations regarding \eqref{Eq:DAPI}.

\renewcommand{\baselinestretch}{1}
\bibliographystyle{IEEEtran}

\bibliography{/Users/jwsimpso/GoogleDrive/JohnSVN/bib/alias,%
/Users/jwsimpso/GoogleDrive/JohnSVN/bib/Main,%
/Users/jwsimpso/GoogleDrive/JohnSVN/bib/JWSP,%
/Users/jwsimpso/GoogleDrive/JohnSVN/bib/New%
}


\end{document}